\documentclass[preprint,a4paper,12pt]{elsarticle} 
\usepackage[]{graphicx}
\usepackage{dsfont}
\newcommand{\bA}{{\bf A}}
\newcommand{\bB}{{\bf B}}
\newcommand{\bC}{{\bf C}}
\newcommand{\bE}{{\bf E}}
\newcommand{\bF}{{\bf F}}
\newcommand{\bG}{{\bf G}}

\newcommand{\bH}{{\bf H}}
\newcommand{\bg}{{\bf g}}\newcommand{\bJ}{{\bf J}}

\newcommand{\bI}{{\bf I}}

\newcommand{\bS}{{\bf S}}
\newcommand{\bX}{{\bf X}}

\newcommand{\ba}{{\bf a}}
\newcommand{\bb}{{\bf b}}
\newcommand{\bc}{{\bf c}}

\newcommand{\tT}{\mathcal{T}}

\journal{Journal of Computational and Applied Mathematics}
\usepackage{hyperref}

\begin{document}
\begin{frontmatter}

\title{Numerical CP Decomposition of Some Difficult Tensors}

\author{Petr Tichavsk\'{y}$^{a,*}$, Anh Huy Phan$^b$ and Andrzej
Cichocki$^b$}

\address{$^{a}$Institute of Information Theory and Automation, Prague 182 08, Czech Republic\\$^{b}$ Brain Science Institute, RIKEN, Wakoshi,
Japan.}

\begin{abstract}
In this paper, a numerical method is proposed for canonical polyadic (CP) decomposition of small size tensors. The focus is primarily on decomposition
of tensors that correspond to small matrix multiplications. Here, rank of the tensors is equal to the smallest number of scalar multiplications that
are necessary to accomplish the matrix multiplication. The proposed method is based on a constrained Levenberg-Marquardt optimization.
Numerical results indicate the rank and border ranks of tensors that correspond to multiplication of matrices of the size $2\times 3$ and $3\times 2$, $3\times 3$ and $3\times 2$,
$3\times 3$ and $3\times 3$, and $3\times 4$ and $4\times 3$. The ranks are 11, 15, 23 and 29, respectively.
In particular, a novel algorithm for multiplying the matrices of the sizes $3\times 3$ and $3\times 2$ with 15 multiplications is presented.
\end{abstract}



\end{frontmatter}

\section{Introduction}

The problem of determining the complexity of matrix multiplication became a well studied topic since the discovery of the Strassen's algorithm \cite{Strassen1}.
The Strassen's algorithm allows multiplying $2\times 2$ matrices using seven multiplications. A consequence of this algorithm is that
$n\times n$ matrices can be multiplied by performing of the order $n^{2.81}$ operations. More recent advances have brought the number of operations needed even closer to the $n^2$ operations.
The current record is $O(n^{2.373})$ operations due to Williams \cite{Williams}.

The problem of the matrix multiplication can be rephrased as a problem of decomposing a particular tensor according to its rank \cite{Landsberg}. The tensor rank
is equal to the lowest number of the scalar multiplications needed to compute the matrix product.
The focus of this paper is not on improving the above asymptotic results but on numerical decomposition of tensors that correspond to
multiplication of small matrices and determining their rank \cite{small}. Although the problem is quite old, only partial results are known so far.

The matrix multiplication tensor for the $2\times 2$ matrices is already completely clear \cite{Winograd}.  Its rank is 7 and its border rank is 7 as well.
The border rank is the lowest rank of tensors that approximate the given tensor. For the $3\times 3$ case, an algorithm for computing the product with 23 scalar multiplications
was found by Laderman \cite{Laderman}.
It means that the rank is at most 23. For multiplying two $4\times 4$ matrices, one can use twice the Strassen's algorithm, and therefore the rank is at most 49.
Multiplication of $5\times 5$ matrices was studied by Makarov \cite{Makarov} with the result of 100 multiplications (rank 100).

In this paper we present a numerical decomposition of the matrix multiplication tensors. For now, we are not able to improve the known results
of Strassen, Laderman and Makarov, we rather show a method of the decomposition with these ranks and numerical results indicating that further improvements
are probably not possible. Moreover, the numerical methods allow to guess the border rank of the tensors.
As a new result, we have derived a novel algorithm for multiplying two matrices of the size $3\times 3$ and $3\times 2$
through 15 multiplications.

Traditional numerical tensor decomposition methods include the alternating least squares method (ALS) \cite{ALS}, 
improved ALS through the
enhanced line search (ELS) \cite{ELS}, damped Gauss-Newton method, also known as Levenberg- Marquardt (LM) method \cite{LM}, and different nonlinear optimization methods, e.g. \cite{NLS}.
For decomposition of the multiplication tensors we have developed a special variant of the constrained LM method. Once an exact fit solution is found, we propose a method of
finding another solution such that the factor matrices only contain nulls, ones and minus ones.

The rest of the paper is organized as follows. The tensors of the matrix multiplication are introduced in Section 2. The numerical method of their decomposition is presented in Section 3.
Section 4 presents numerical results and section 5 concludes the paper.

\section{Tensor of Matrix Multiplication}

Consider two matrices $\bE$ and $\bF$ of the sizes $P\times Q$ and $Q\times S$, respectively, and their matrix product
$\bG=\bE\bF$ of the size $P\times S$. The operation of the matrix multiplication can be represented by a tensor $\tT_{PQS}$ of the size
$PQ\times QS\times PS$ which is filled with nulls and ones only, such that
\begin{equation}
\mbox{vec}(\bG)=\tT_{PQS}\times_1\mbox{vec}(\bE^T)^T\times_2\mbox{vec}(\bF^T)^T\label{jedna}
\end{equation}
regardless of the elements values of $\bE$ and $\bF$. Here, $\times_i$ denotes a tensor-matrix multiplication along the dimension $i$, and
$\mbox{vec}$ is an operator that stacks all elements of a matrix or tensor in one long column vector.

Note that the number of ones in the tensor $\tT_{PQS}$ is $PQS$; it is the number of scalar multiplications needed for evaluating
the matrix product by a conventional matrix multiplication algorithm.

The tensor $\tT_{PQS}$ can be obtained by reshaping an order-6 tensor with elements
\begin{equation}
{\mathcal T}^{(PQS)}_{ijk\ell mn}=\delta_{i\ell}\delta_{jm}\delta_{kn}\qquad\mbox{for}\quad i,k=1,\ldots,P; \ell,m=1,\ldots,R; j,n=1,\ldots,S\label{six}
\end{equation}
to the format $PQ\times QS\times PS$.

For example,
\begin{equation}
{\mathcal T}_{222}=\left(\begin{array}{cccc}  1  &   0  &   0 & 0\\
     0 &    0 &    1 & 0\\ 0& 0& 0& 0\\ 0&0 &0 &0\end{array}\left\vert
     \begin{array}{cccc}  0&   0  &   0 & 0\\0& 0& 0& 0\\
     1 &    0 &    0 & 0\\  0&0 &1 &0\end{array}\right\vert\left.
     \begin{array}{cccc}   0  &   1  &   0 & 0\\
     0 &    0 &  0& 1\\ 0& 0& 0& 0\\ 0&0 &0 &0\end{array}\right\vert
     \begin{array}{cccc}  0&   0  &   0 & 0\\0& 0& 0& 0\\
     0  &   1  &   0 & 0\\  0 &    0 &  0& 1\end{array}
     \right)~.\label{222}
\end{equation}
This tensor has the size $4\times 4\times 4$, the vertical lines separate the four frontal slices of the tensor.
\vspace{5mm}

A canonical polyadic decomposition of the tensor $\tT_{PQS}$ is a representation of the tensor as a sum of $R$ rank-one components
\begin{eqnarray*}
\tT_{PQS}=\sum_{r=1}^R \ba_r\circ\bb_r\circ \bc_r
\end{eqnarray*}
where $\{\ba_r\}$, $\{\bb_r\}$, $\{\bc_r\}$ are columns of so called factor matrices $\bA,\bB,\bC$, symbolically $\tT_{PQS}=[[\bA,\bB,\bC]]$.
For example, a CP decomposition of the tensor $\tT_{222}$ in (\ref{222}) corresponding to the Strassen algorithm [2] is ${\mathcal T}_{222}=[[\bA,\bB,\bC]]$ with
\begin{eqnarray*}
\bA&=&\left(\begin{array}{ccccccc}  1& 0& 1& 0& 1& -1& 0\\ 0& 0& 0& 0& 1& 0& 1\\ 0& 1& 0& 0& 0& 1& 0\\ 1& 1& 0& 1& 0& 0& -1
   \end{array}\right)\\
\bB&=&\left(\begin{array}{ccccccc}  1& 1& 0& -1& 0& 1& 0\\ 0& 0& 1& 0& 0& 1& 0\\ 0& 0& 0& 1& 0& 0& 1\\ 1& 0& -1& 0& 1& 0& 1
   \end{array}\right)\\
\bC&=&\left(\begin{array}{ccccccc}  1& 0& 0& 1& -1& 0& 1\\  0& 1& 0& 1& 0& 0 &0\\ 0& 0& 1& 0& 1& 0& 0\\ 1& -1& 1& 0& 0& 1 &0 \end{array}\right)
\end{eqnarray*}

The multiplication tensors have the following properties:
\begin{enumerate}
\item Ranks of these tensors exceed the tensors' dimensions.
\item The CP decompositions are not unique.
\item The border ranks of the tensors might be strictly lower than their true ranks.
\item Tensors $\tT_{NNN}$ are invariant with respect to some permutations of indices. Using the matlab notation
we can write
${\mathcal T}_{NNN}=\mbox{permute}({\mathcal T}_{NNN},[2,3,1])=\mbox{permute}({\mathcal T}_{NNN},[3,1,2])$
\item Tensors $\tT_{PQS}$ are invariant with respect to certain tensor-matrix multiplications [1].
\end{enumerate}
Let us explain the last item in  more details.
Since it holds $\bG=\bE\bF=(\bE\bX)(\bX^{-1}\bF)$ for any invertible matrix $\bX$, we have
$$ 
\mbox{vec}(\bG)=\tT_{PQS}\times_1\mbox{vec}(\bE^T)^T\times_2\mbox{vec}(\bF^T)^T=\tT_{PQS}\times_1\mbox{vec}(\bX^T\bE^T)^T\times_2\mbox{vec}(\bF^T\bX^{-T})^T~.
$$ 
The multiplication with $\bX$ and $\bX^{-1}$ can be absorbed into $\tT_{PQS}$, because
\begin{eqnarray*}
\mbox{vec}(\bX^T\bE^T)&=&(\bI\otimes\bX^T)\mbox{vec}(\bE^T)\\
\mbox{vec}(\bF^T\bX^{-T})&=&(\bX^{-1}\otimes\bI)\mbox{vec}(\bF^T)
\end{eqnarray*}
where $\bI$ is identity matrix of an appropriate size. Therefore
$$\tT_{PQS}=\tT_{PQS}\times_1 \bS_1(\bX)\times_2 \bS_2(\bX)$$ where $\bS_1(\bX)=\bI\otimes \bX^T$  and $\bS_2(\bX)=
\bX^{-1}\otimes \bI$.

\section{Numerical CP Decomposition}

For numerical CP decomposition of the multiplication tensors we propose a three--step procedure.
\begin{enumerate}
\item Finding an ``exact fit'' solution, if it exists. \item  Finding another solution which would be as much sparse as possible. \item  Finding a solution with factor matrices containing only integer (or rational) entries.
\end{enumerate} \vspace{5mm}

{\bf Step 1: Finding an Exact Fit} \vspace{5mm}

We seek a vector of parameters $\theta_R=[\mbox{vec}({\bf A})^T,\mbox{vec}({\bf B})^T,\mbox{vec}({\bf C})^T]^T$
of the size $3N^2R\times 1$ that minimizes the cost function $\varphi(\theta)=\|{\mathcal T}_N-\hat{\mathcal T}(\theta)\|_F^2$ in the compact set
$\{\theta\in\mathds{R}^{3N^2R}; \|\theta\|^2=c\}$, where $c$ is a suitable constant.

The ordinary (unconstrained) LM algorithm updates $\theta$
as
\begin{equation}
\theta \leftarrow \theta - (\bH+\mu \bI)^{-1}{\bf g}
\end{equation}
where
\begin{equation}
\bH=\bJ^T\bJ,\qquad \bJ=\frac{\partial\mbox{vec}(\hat{\tT}(\theta))}{\partial\theta},\qquad \bg=\bJ^T\mbox{vec}(\tT_N-\hat{\tT}(\theta))
\end{equation}
and $\mu$ is a damping parameter, which is sequentially updated according to a rule described in \cite{Madsen}. Closed-form expressions for the
Hessian $\bH$ and gradient $\bg$ can be found in \cite{LM}.

Optimization constrained to the ball is performed by minimizing the cost function in the tangent plane $\{\theta; (\theta-\theta_0)^T\theta_0=0\}$ first,
where $\theta_0$ is the latest available estimate of $\theta$ which obeys the constraint.
Let $\theta_1^\prime$ be the minimizer in the tangent plane.
Then, $\theta_1^\prime$ is projected on the ball by an appropriate scale change,
$\theta_1=\theta_1^\prime\sqrt{c}/\|\theta_1^\prime\|$.\vspace*{2mm}

Towards computing $\theta_1^\prime$,
let the following second-order approximation of the cost function be minimized,
\begin{equation}
\varphi(\theta)\approx\varphi(\theta_0)+{\bf g}^T(\theta-\theta_0)+\frac{1}{2}(\theta-\theta_0)^T\bH(\theta-\theta_0)
\end{equation}
under the linear constraint $(\theta-\theta_0)^T\theta_0=0$. We use the method of Lagrange multiplier to get
\begin{equation}
\theta_1^\prime=\theta_0-\bH^{-1}{\bf g}+\frac{\theta_0^T\bH^{-1}{\bf g}}{\|\theta_0\|^2}\,\bH^{-1}\theta_0~.\label{copt}
\end{equation}
Instead of using (4) directly, we propose replacing $\bH^{-1}$ by $(\bH+\mu \bI)^{-1}$ as in the LM method.
\vspace{2mm}

We need to do multiple random initializations to get close to the global minimum of the cost function; in the optimum case it is the
exact fit solution, i.e. with $\varphi(\theta)=0$. The method works well for small matrices. For example, for decomposition
of the $\tT_{333}$ and constraint $c=150$ we need only a few random trials to obtain an exact fit solution. On the other hand, for tensor
$\tT_{444}$ the false local minima are so numerous that it is almost impossible to get an exact fit decomposition
when the algorithm is started from random initial conditions.
\vspace{5mm}

{\bf Step 2: Finding a Sparse Solution}\vspace{5mm}

For simplicity, we describe a method of finding a sparse CP decomposition in the case of tensors $\tT_{NNN}$.
Let $\tT_{NNN}=[[\bA,\bB,\bC]]$ be an exact CP decomposition with certain $\bA,\bB,\bC$. We have
\begin{equation}
\tT_{NNN}=[[\bS_1(\bX)\bA,\bS_2(\bX)\bB,\bC]]
\end{equation}
where $\bX$ is an arbitrary invertible matrix of size $N\times N$.

First, we seek a matrix $\bX$ of determinant 1 such that $\|\bS_1(\bX)\bA\|_1+\|\bS_2(\bX)\bB\|_1$ is minimized, and update $\bA,\bB$
as $\bA\leftarrow \bS_1(\bX)\bA$, $\bB\leftarrow \bS_2(\bX)\bB$. We use the Nelder-Mead algorithm for the minimization.

Second, we seek another $\bX$ such that
$\|\bS_1(\bX)\bB\|_1+\|\bS_2(\bX)\bC\|_1$ is minimized,
and update $\bB$ and $\bC$.

Third, we seek another $\bX$ such that $\|\bS_1(\bX)\bC\|_1+\|\bS_2(\bX)\bA\|_1$ is minimized,
and update $\bC$ and $\bA$.

The sequence of three partial optimizations is repeated until convergence is obtained.
\vspace{5mm}

As a result, we obtain $\tT_{NNN}=[[\bA,\bB,\bC]]$ where many elements of $\bA,\bB,\bC$ are nulls.
\vspace{5mm}

{\bf Step 3: Finding a Rational Solution}\vspace{5mm}

We continue to modify the exact fit solution obtained in the previous step by constraining some other elements of $\theta_R$ to be 1 or -1. We do this
by sequentially increasing the number of elements of $\theta_R$ to be in the set $\{0,1,-1\}$.
In each step, the function
$\varphi(\theta_R)$ is minimized, starting from the latest available solution, with another free element of $\theta_R$ changed and fixed to 1 or -1. If an exact fit cannot be achieved, another free element is tried instead.
At the very end, it might happen that none of the free elements of $\theta_R$ can be set to 1 or -1. In that case, we suggest to try the values 2 or -2 or higher.
Some other elements of $\theta_R$ may become $1/2$ or $-1/2$.

\section{Experiments}

\subsection{Estimating the tensor rank}
For the multiplication tensor $\tT_{333}$, holds $\min \varphi(\theta_{23})=0$ under the constraint $\|\theta_{23}\|^2=150$. The exact fit can be obtained quite quickly.
For a rank--22 approximation of $\tT_{333}$, even with a more relaxed constraint $\|\theta_{22}\|^2=594$, the lowest possible value of the fit that we were able to find was $\min \varphi(\theta_{22})=6.766\cdot 10^{-5}>0$~. These observations indicate that the rank of the tensor ${\mathcal T}_{333}$ is 23. Similarly, if we attempt to decompose the tensor
to 22 and 21 terms, $\min \varphi(\theta)$ converges to zero for $c\rightarrow \infty$. However, for decomposition to 20 terms, $\min \varphi(\theta)$ does not converge to zero.
Therefore we make the conjecture that the border rank of the tensor is 21.

A more complete table of numerical results obtained by the above described procedure is as follows.
\begin{table}
\begin{center}
TABLE I\\
{\small Upper bounds for ranks and border ranks of multiplication tensors}
\begin{tabular}{|c|c|c|c|c|c}\hline acronym & matrix sizes & \# of 1's & rank & border rank
 \\ \hline
$222$  & $2\times2$,\quad $2\times 2$ & 8 & 7 & 7 \\
$232$  & $2\times3$,\quad $3\times 2$ & 12 & 11 & 10 \\
$322$  & $3\times2$,\quad $2\times 2$ & 12 & 11 & 10 \\
$332$  & $3\times3$,\quad $3\times 2$ & 18 & 15 & 14 \\
$333$  & $3\times3$,\quad $3\times 3$ & 27 & 23 & 21 \\
$343$  & $3\times4$,\quad $4\times 3$ & 36 & 29 & 28 \\
$443$  & $4\times4$,\quad $4\times 3$ & 48 & 40 & 39 \\
$444$  & $4\times4$,\quad $4\times 4$ & 64 & 49 & 49 \\
\hline
\end{tabular}
\end{center}
\end{table}
The table shows rank of the exact-fit solutions of the CP decomposition of the tensors obtained by the constrained optimization.
The numerical border rank was determined as the minimimum rank for which min$\varphi(\theta)$ constrained by $\|\theta\|=c$
converges to zero for $c\rightarrow \infty$. It is not a mathematical proof that the border rank has the displayed values,
but an empirical observation. The true border ranks can be theoretically smaller.

The results in the table are rather discouraging for multiplication of the matrices $2\times3$ with $3\times 2$ and
$3\times2$ with $2\times 2$. Our experiments indicate that the necessary number of multiplications is 11 in these two cases.
Corresponding algorithm can be obtained
by applying the Strassen algorithm to the $2\times2$ blocks. It is not interesting from the computational point of view.

The only novel results are obtained for the cases $3\times3$ with $3\times 2$ and
$3\times 4$ with $4\times 3$. The former case is studied in the next subsection in more details.

\includegraphics[width=\linewidth,height=0.5\linewidth]{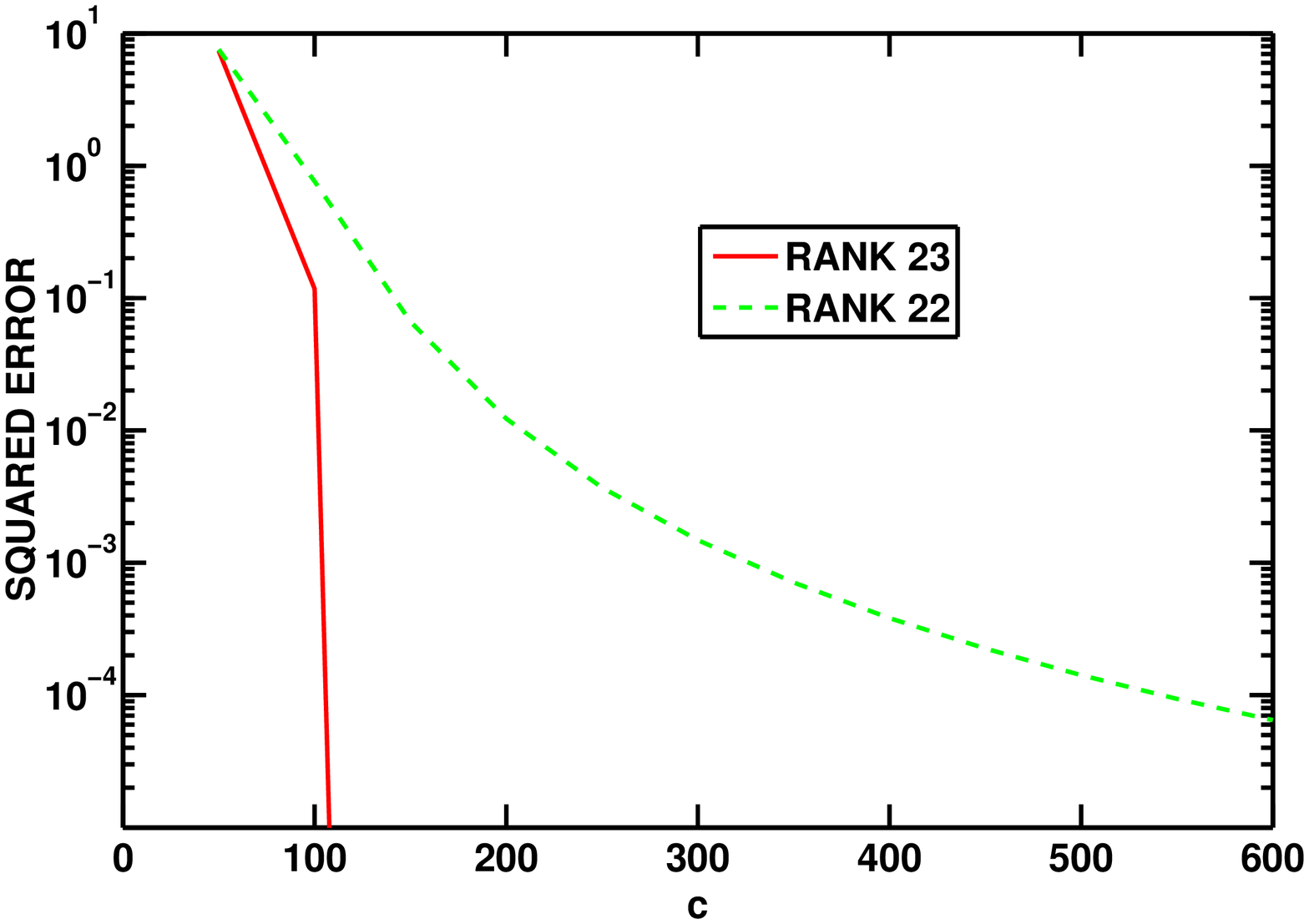}
\vspace{5mm}
Fig. 1: Reconstruction error as a function of parameter $c$ in the constraint $\|\theta_R\|^2=c$ for $N=3$, $R=22$ and $R=23$.

\includegraphics[width=\linewidth,height=0.5\linewidth]{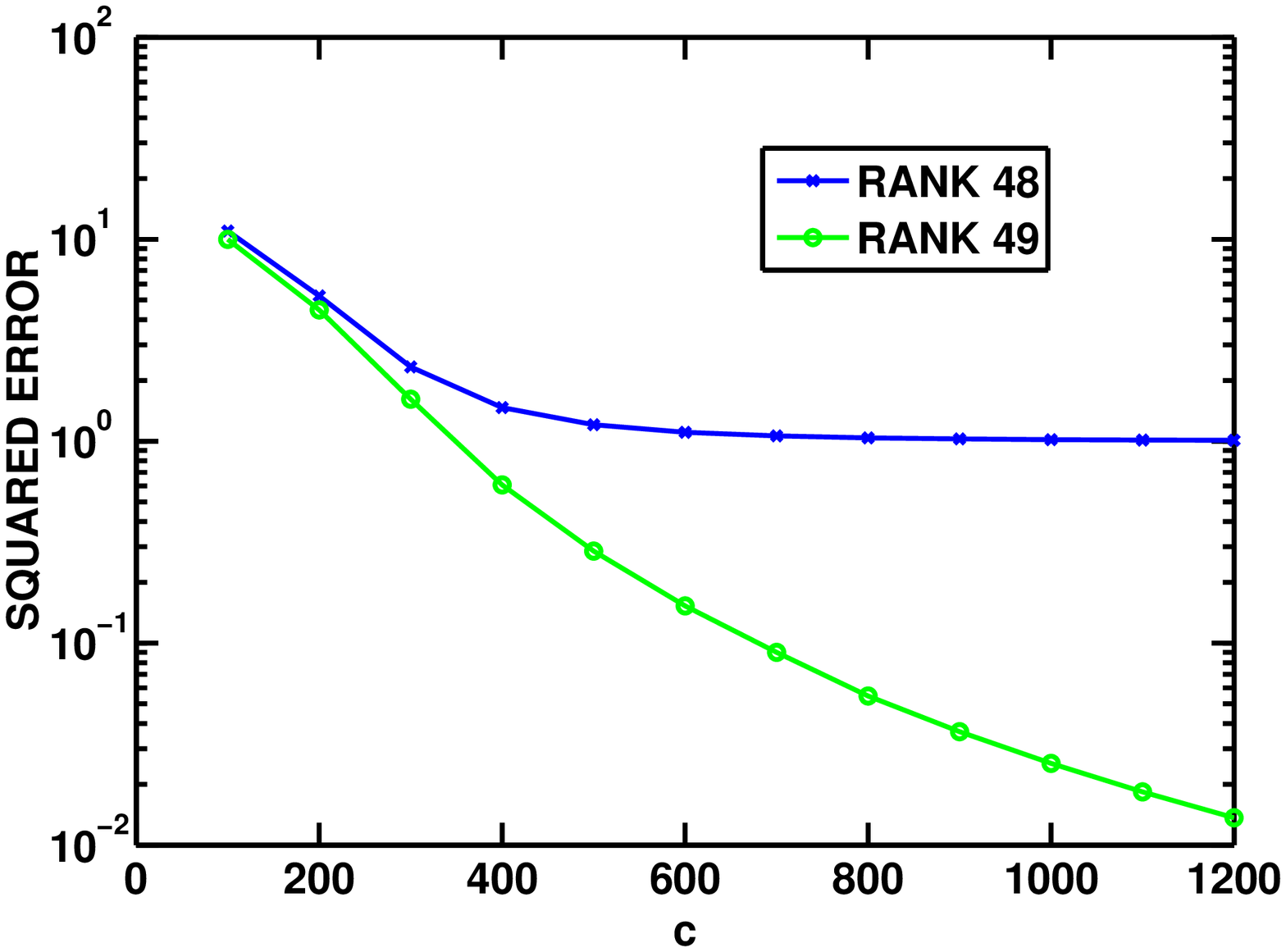}
\vspace{5mm}
Fig. 2: Reconstruction error as a function of parameter $c$ in the constraint $\|\theta_R\|^2=c$ for $N=4$, $R=48$ and $R=49$.

\subsection{Multiplication of Matrices $3\times 3$ and $3\times 2$}

Consider multiplication of matrices $3\times 3$ and $3\times 2$ of the form
\begin{equation}
\left(\begin{array}{ccc} e_{11} & e_{12} & e_{13}\\e_{21} & e_{22} & e_{23}\\ e_{31} & e_{32} & e_{33}\end{array}\right)
\left(\begin{array}{cc} f_{11} & f_{12} \\ f_{21} & f_{22}\\ f_{31} & f_{32}\end{array}\right)=
\left(\begin{array}{cc} g_{11} & g_{12} \\ g_{21} & g_{22} \\ g_{31} & g_{32}\end{array}\right)\label{product}
\end{equation}
Standard algorithm for computing $g_{11},\ldots,g_{32}$ from $\{e_{ij}\}$ and $\{f_{ij}\}$
requires 18 scalar multiplications.
We show that the computation can be accomplished through 15 scalar multiplications only.\\
The tensor representing the multiplication has dimension $9\times 6\times 6$.
A CP decomposition of this tensor obtained by applying the proposed algorithm is ${\mathcal T}_{332}=[[A,B,C]]$, where
\begin{eqnarray*}
A&=&\left(\begin{array}{ccccccccccccccc}    
    -1  &   0 &    0 &    0 &    0 &   -1  &   0  &   0  &   0  &   0  &   0 &    0  &   0  &   0 &    0\\
    -1  &   0 &    0 &   -1 &    0 &    0  &   0  &   1   &  0  &   0  &   1 &    0  &   0  &   0 &    0\\
     0  &  -1 &    0 &    0 &   -1 &    1  &   0  &   0  &   0  &   0  &   0 &    0  &   0  &   0 &    1\\
     0  &   0 &    1 &    1 &    0 &   -1  &   0  &   0  &   0  &   0  &   0 &    1  &   0  &   0 &    0\\
     0  &   0 &    0 &    0 &    0 &    0  &   0  &   0  &   0  &   0  &   1  &   1  &   0  &   0 &    0\\
     0  &   0 &    0 &    0 &    0 &    0  &   0  &   0  &   1  &  -1  &   1  &   0  &  -1  &   0 &    0\\
     1  &   1 &    0 &    0 &    0 &    0  &   1  &   0  &   0  &  -1  &   0  &   0  &   0  &   0 &    0\\
     0  &   0 &    0 &    0 &    1 &    0  &   0  &   0  &   0  &   0  &   0  &  -1  &  -1  &  -1 &    0\\
     0  &   0 &    0 &    0 &    1 &    0  &   0  &   0  &   0  &   1   &  0  &   0  &   0  &   0 &    0\end{array}\right)\\
B&=&\left(\begin{array}{ccccccccccccccc}  
     1   &  1   &  0   &  1   &  0 &   -1  &   0  &  -1  &   0  &   0  &   0  &   0 &    0  &   0  &   1\\
    -1   &  0   & -1   & -1   &  0 &    0  &   1  &   1  &   0  &   0  &   0   &  0 &    0  &   0  &   0\\
     0   &  0   &  0   &  0   &  0 &    0  &   0  &   1  &   1  &   0  &   1  &   0 &   -1  &  -1  &   0\\
     0   &  0   &  1   &  1   &  0 &    0  &   0  &  -1  &   0  &   0  &   0  &   1 &    0  &   1  &   0\\
     0   &  0   &  0   &  0   &  1 &    0  &   0  &   0  &  -1  &   0  &   0  &   0 &    1   &  1  &   1\\
     0   &  1   &  0   &  0   & -1 &    0  &   1  &   0  &   0  &   1  &   0  &   0 &   -1  &  -1  &   0\end{array}\right)\\
C&=&\left(\begin{array}{ccccccccccccccc} 
     0  &   0 &    1 &   -1 &    0 &    1  &   0 &    1  &   0 &    0   &  0  &   0  &   0  &   0   &  1\\
     0  &   0 &   -1 &    1 &    0 &    0  &   0 &   -1  &  -1 &    0   &  1  &   0  &   0  &   0   &  0\\
     0  &   1 &    0 &    0 &    1 &    0  &   0 &    0  &  -1 &    1   &  0  &   0  &   1  &   0   &  1\\
     1  &  -1 &    1 &   -1 &    0 &    1  &   1 &    0  &   0 &    0   &  0   &  0  &   0  &   0   &  0\\
     0  &   0 &   -1 &    0 &    0  &   0  &   0 &    0  &  -1 &    0   &  0  &   1  &   1  &  -1   &  0\\
     0  &   0 &    0 &    0 &    0 &    0  &   1 &    0  &  -1 &    1   &  0  &   0   &  1  &  -1   &  0\end{array}\right)
\end{eqnarray*}
The 15 scalar multiplications are
\begin{eqnarray*}
m_1 &=& -(f_{11} - f_{12})(e_{11} + e_{12} - e_{31})\\
m_2 &=&       -(f_{11} + f_{32})(e_{13} - e_{31})\\
m_3 &=&               -e_{21}(f_{12} - f_{22})\\
m_4 &=& -(e_{12} - e_{21})(f_{11} - f_{12} + f_{22})\\
m_5 &=&  (f_{31} - f_{32})(e_{32} - e_{13} + e_{33})\\
m_6 &=&          f_{11}(e_{11} - e_{13} + e_{21})\\
m_7 &=&                e_{31}(f_{12} + f_{32})\\
m_8 &=&   -e_{12}(f_{11} - f_{12} - f_{21} + f_{22})\\
m_9 &=&                e_{23}(f_{21} - f_{31})\\
m_{10} &=&         -f_{32}(e_{23} + e_{31} - e_{33})\\
m_{11} &=&          f_{21}(e_{12} + e_{22} + e_{23})\\
m_{12} &=&          f_{22}(e_{21} + e_{22} - e_{32})\\
m_{13} &=&  (e_{23} + e_{32})(f_{21} - f_{31} + f_{32})\\
m_{14} &=&    e_{32}(f_{21} - f_{22} - f_{31} + f_{32})\\
m_{15} &=&                e_{13}(f_{11} + f_{31})
\end{eqnarray*}
Having computed these products, the elements $g_{ij}$ in (\ref{product}) can be written as
\begin{eqnarray*}
\left(\begin{array}{cc} g_{11} & g_{12} \\ g_{21} & g_{22} \\ g_{31} & g_{32}\end{array}\right)=\left(\begin{array}{cc} m_3-m_4+m_6+m_8+m_{12} & m_1-m_2+m_3-m_4+m_6+m_7 \\ -m_3+m_4-m_8-m_9+m_{11} & -m_3-m_9+m_{12}+m_{13}-m_{14} \\ m_2+m_5-m_9+m_{10}+m_{13}+m_{15} & m_7 -m_9+m_{10}+m_{13}-m_{14}\end{array}\right)
\end{eqnarray*}

\section{Conclusions}

 The constrained LM algorithm may serve for decomposition of difficult tensors that have the border rank lower than the true rank
and when uniqueness is not required. Numerical decomposition of tensors larger than $\tT_{333}$, e.g. $\tT_{444}$, is still a challenging task.
We have provided a decomposition of $\tT_{332}$ tensor to 15 rank-one terms, i.e. showed that product of the matrices $3\times 3$ and $3\times 2$ can be computed through 15 scalar multiplications.
Matlab codes of the proposed algorithms are available on the web page of the first author.

\section*{Acknowledgements}
The work of Petr Tichavsk\'y was
supported by Czech Science Foundation through project 14-13713S.


\begin{thebibliography}{99}

\bibitem{Strassen1}  V. Strassen, Gaussian elimination is not optimal, Numer. Math. 13 (1969) 354--356.

\bibitem{Williams} V.V.~Williams, Multiplying Matrices Faster Than Coppersmith-Winograd, in Proc. of
 the 44th Symposium on Theory of Computing, STOC '12, New York, NY, USA, (2012) 887-–898.

\bibitem{Landsberg} J. M. Landsberg, Tensors: Geometry and Applications, AMS 2012.

\bibitem{small} C.-E. Drevet, M.N. Islam and E. Schost, Optimization techniques for small matrix multiplications, Theoretical Computer Science 412 (2011), 2219-–2236.

\bibitem{Winograd} S. Winograd, On Multiplication of $2\times 2$ Matrices, Linear Algebra and Appl. 4 (1971) 381-–388.

\bibitem{Laderman}  J.D. Laderman, A noncommutative algorithm for multiplying $3\times 3$ matrices using 23 multiplications, Bul. Amer. Math. Soc. 82 (1976) 126--128.

\bibitem{Makarov}  O.M. Makarov, A noncommutative algorithm for multiplying $5\times 5-$ matrices using one hundred multiplications, U.S.S.R. Comput. Maths. Math. Phys. 27 (1987) 311--315.

\bibitem{ALS} P. Comon, X. Luciani and A.~L.~F.~de~Almeida, Tensor decompositions, alternating least
squares and other tales, Chemometrics 23 (2009) 393-405.

\bibitem{ELS} M. Rajih, P. Comon, and R.~A.~Harshman,
Enhanced line search: A novel method to accelerate PARAFAC, SIAM Journal on Matrix Analysis Appl. 30 (2008) 1148--1171.


\bibitem{LM} A.H. Phan, P. Tichavsk\'y and A. Cichocki,
Low Complexity Damped Gauss-Newton Algorithms for Parallel Factor
Analysis, SIAM J. Matrix Anal. and Appl. 34 (2013) 126--147.


\bibitem{NLS} L. Sorber, M. Van Barel, and L. De Lathauwer, Structured data fusion,
Tech. Rep., ESAT-SISTA, Internal Report 13-177, 2013.

\bibitem{Madsen} K. Madsen, H. B. Nielsen, O. Tingleff, Methods for nonlinear
least squares problems, second ed., Department of Mathematical
Modelling, Technical University of Denmark, Lyngby,
Denmark, 2004.


\end{thebibliography}
\end{document}